%%%%%%%%%%%%%%%%%%%%%%%%%%%%%%%%%%%%%%%%%%%%%%%%%%%%%%%%%%%%%%%%%%%%%%%%%%%%%%%%%%%%%%%
%%%%%%%%%%%%%%%%%%%%%%%%%%%%%%%%%%%%%%%%%%%%%%%%%%%%%%%%%%%%%%%%%%%%%%%%%%%%%%%%%%%%%%%
%%%%%% The c-function for non-compactly causal symmetric spaces %%%%%%%%%%%%%%%%%%%%%%% 
%%%%%%%%%%%%%%%%%%%%%%%%%%%%%%%%%%%%%%%%%%%%%%%%%%%%%%%%%%%%%%%%%%%%%%%%%%%%%%%%%%%%%%%
%%%%%%%%%%%%%%%%%%%%%%%%%%%%%%%%%%%%%%%%%%%%%%%%%%%%%%%%%%%%%%%%%%%%%%%%%%%%%%%%%%%%%%%

%VERSION OF august 24, 1995
\input amssym.def
\input amssym.tex
% Achtung: Bei Verwendung von jltmac werden die Fontfamilien masm und
% msbm mehrfach aufgerufen. Das kann zu Kapazitaetsschwierigkeiten
% fuehren. KH 27.2.95

%\def\lit{\sectionheadline{\bf References}
%\frenchspacing
%\entries\par}

\def\item#1{\vskip1.3pt\hang\textindent {\rm #1}}% THIS REPLACES KNUTH'S DEF'N

                                  % THIS REPLACES KNUTH'S DEF'N

\tolerance=300
\pretolerance=200
\hfuzz=1pt
\vfuzz=1pt

% Print out with -x300 -y50

% Offsetwerte fuer Ausdrucke in Erlangen, hoffset um 0.6 in groesser machen
\hoffset=0.6in
\voffset=0.8in

%\baselineskip                 
\hsize=5.8 true in 

%%In Deutschland \vsize=9.2 %%

\vsize=8.5 true in
%\baselineskip
\parindent=25pt
\mathsurround=1pt
\parskip=1pt plus .25pt minus .25pt
\normallineskiplimit=.99pt

\countdef\revised=100
\mathchardef\emptyset="001F % THIS REPLACES KNUTH'S DEFINITION
\chardef\ss="19
\def\3{\ss}
\def\anf{$\lower1.2ex\hbox{"}$}
\def\frac#1#2{{#1 \over #2}}
\def\>{>\!\!>}
\def\<{<\!\!<}

\def\ssarr{\hbox to 30pt{\rightarrowfill}}
\def\sarr{\hbox to 40pt{\rightarrowfill}}
\def\arr{\hbox to 60pt{\rightarrowfill}}
\def\larr{\hbox to 60pt{\leftarrowfill}}
\def\Arr{\hbox to 80pt{\rightarrowfill}}

{}

\def\ad{\mathop{\rm ad}\nolimits}

\def\Ad{\mathop{\rm Ad}\nolimits}

\def\cone{\mathop{\rm cone}\nolimits}

\def\ev{\mathop{\rm ev}\nolimits}

%\def\Fa{\mathop{\rm Fa}\nolimits}

%
%
 % USED FOR IDENTITY FUNCTION

% USED FOR IMAGINARY PART OF COMPLEX NUMBERS

\def\Inn{\mathop{\rm Inn}\nolimits}
\def\Int{\mathop{\rm int}\nolimits}

\def\Re{\mathop{\rm Re}\nolimits}% USED FOR REAL PART OF COMPLEX NUMBERS

\def\Sl{\mathop{\rm Sl}\nolimits}
\def\SO{\mathop{\rm SO}\nolimits}
\def\span{\mathop{\rm span}\nolimits}

\def\Spec{\mathop{\rm Spec}\nolimits}

% USED FOR TRACE OF MATRIX
%\def\trdeg{\mathop{\rm trdeg}\nolimits} 

%\def\Vequiv{\mathrel{_V\!\equiv}}

\def\0{{\bf 0}}
\def\1{{\bf 1}}

\def\a{{\frak a}}

\def\b{{\frak b}}
\def\cc{{\frak c}}

\def\g{{\frak g}}

\def\h{{\frak h}}

\def\k{{\frak k}}

\def\n{{\frak n}}

\def\p{{\frak p}}
\def\q{{\frak q}}

\def\z{{\frak z}}

\def\C{{\Bbb C}}

\def\N{{\Bbb N}}

\def\R{{\Bbb R}}

\def\:{\colon}  %8.5.92
\def\.{{\cdot}}
\def\|{\Vert}
\def\bsk{\bigskip}

\def\giantskip{\vskip2\bigskipamount}
\def\gsk{\giantskip}
\def \la {\langle}
\def\msk{\medskip}
\def \ra {\rangle}
\def \res {\!\mid\!\!}

\def\ssk{\smallskip}

\def\bbr{\bigbreak}
\def\giantbreak{\par \ifdim\lastskip<2\bigskipamount \removelastskip
         \penalty-400 \giantskip\fi}

\def\nin{\noindent}
\def\cen{\centerline}
\def\pagebreak{\vskip 0pt plus 0.0001fil\break}
\def\linebreak{\break}

\def\eps{\varepsilon}
\def\epsilon{\varepsilon}
\def\eset{\emptyset}

\def\nin{\noindent}
\def\oline{\overline}

\def\pder#1,#2,#3 { {\partial #1 \over \partial #2}(#3)}
\def\pde#1,#2 { {\partial #1 \over \partial #2}}
\def\phi{\varphi}

% Besser \ltimes und \rtimes aus dem AMS-Symbols verwenden. 
%\def\sdir#1{\hbox{$\mathrel\times{\hskip -4.6pt 
%            {\vrule height 4.7 pt depth .5 pt}}\hskip 2pt_{#1}$}}

\def\subeq{\subseteq}
\def\supeq{\supseteq}

\def\Rarrow{\Rightarrow}

\def\tilde{\widetilde}

\font\eightrm=cmr8

% SANS SERIF 10 POINT
 %SANS SERIF 10 POINT ITALIC

\font\smc=cmcsc10
%\font\smc8=cmcsc8 
 %SLANTED TYPEWRITER 10 POINT
 %BOLD FACE MATH SYMBOLS 10 POINT
 %DUNHILL STYLE 10 POINT
 %SAN SERIF BOLD EXTENDED 10 POINT
 %USED FOR TITLES
 %USED FOR TITLES
\font\bfone=cmbx10 scaled\magstep1 %BOLDFACE AT MAGSTEP 1
\font\bftwo=cmbx10 scaled\magstep2 %BOLDFACE AT MAGSTEP 2
 %BOLDFACE AT MAGSTEP 3

\def\qed{{\unskip\nobreak\hfil\penalty50\hskip .001pt \hbox{}\nobreak\hfil
          \vrule height 1.2ex width 1.1ex depth -.1ex
           \parfillskip=0pt\finalhyphendemerits=0\medbreak}\rm}
%This is the end-of-proof sign. 
%Not to be used in display mode. 
%If you want to conclude a proof 
%at the end of a line is display mode use 

%BUT OMIT $$---the macro will write that

\def\Lemma #1. {\bigbreak\vskip-\parskip\noindent{\bf Lemma #1.}\quad\it}

\def\Sublemma #1. {\bigbreak\vskip-\parskip\noindent{\bf Sublemma #1.}\quad\it}

\def\Proposition #1. {\bigbreak\vskip-\parskip\noindent{\bf Proposition #1.}
\quad\it}

\def\Corollary #1. {\bigbreak\vskip-\parskip\nin{\bf Corollary #1.}
\quad\it}

\def\Theorem #1. {\bigbreak\vskip-\parskip\noindent{\bf Theorem #1.}
\quad\it}

\def\Definition #1. {\rm\bigbreak\vskip-\parskip\noindent{\bf Definition #1.}
\quad}

\def\Remark #1. {\rm\bigbreak\vskip-\parskip\noindent{\bf Remark #1.}\quad}

\def\Example #1. {\rm\bigbreak\vskip-\parskip\noindent{\bf Example #1.}\quad}

\def\Problems #1. {\bigbreak\vskip-\parskip\noindent{\bf Problems #1.}\quad}
\def\Problem #1. {\bigbreak\vskip-\parskip\noindent{\bf Problems #1.}\quad}

\def\Conjecture #1. {\bigbreak\vskip-\parskip\noindent{\bf Conjecture #1.}\quad}

\def\Proof#1.{\rm\par\ifdim\lastskip<\bigskipamount\removelastskip\fi\smallskip
            \noindent {\bf Proof.}\quad}

\def\Axiom #1. {\bigbreak\vskip-\parskip\noindent{\bf Axiom #1.}\quad\it}

\def\Satz #1. {\bigbreak\vskip-\parskip\noindent{\bf Satz #1.}\quad\it}

\def\Korollar #1. {\bbr\vskip-\parskip\nin{\bf Korollar #1.} \quad\it}

\def\Bemerkung #1. {\rm\bigbreak\vskip-\parskip\noindent{\bf Bemerkung #1.}
\quad}

\def\Beispiel #1. {\rm\bigbreak\vskip-\parskip\noindent{\bf Beispiel #1.}\quad}
\def\Aufgabe #1. {\rm\bigbreak\vskip-\parskip\noindent{\bf Aufgabe #1.}\quad}

\def\Beweis#1. {\rm\par\ifdim\lastskip<\bigskipamount\removelastskip\fi
           \smallskip\noindent {\bf Beweis.}\quad}

\nopagenumbers

\def\date{\ifcase\month\or January\or February \or March\or April\or May
\or June\or July\or August\or September\or October\or November
\or December\fi\space\number\day, \number\year}

\def\title{Title ??}
\def\author{Author ??}

\def\thanks#1{\footnote*{\eightrm#1}}

\def\rightheadline{\hfil{\eightrm\title}\hfil\tenbf\folio}
\def\leftheadline{\tenbf\folio\hfil{\eightrm\author}\hfil}
\headline={\vbox{\line{\ifodd\pageno\rightheadline\else\leftheadline\fi}}}

\def\firstheadline{}
\def\firstfootline{\cen{\rm\folio}}

\def\seite #1 {\pageno #1
               \headline={\ifnum\pageno=#1 \firstheadline
               \else\ifodd\pageno\rightheadline\else\leftheadline\fi\fi}
               \footline={\ifnum\pageno=#1 \firstfootline\else{}\fi}}

%%%THIS IS THE MACRO LEFTSPACE.TEX %%%TO THD VIA WAFRUPP
\newdimen\dimenone
 \def\checkleftspace#1#2#3#4{%DIESER MACRO STAMMT VON APPELT
 \dimenone=\pagetotal%#1=Skip vorher,#2=Font,#3=Text,#4=Skip nachher  
 \advance\dimenone by -\pageshrink   %testen ob Titel noch mit Gewalt auf Seite 
                                                                          %geht
 \ifdim\dimenone>\pagegoal          %nacha tua nix-- gewoehnliche Outputroutine 
   \else\dimenone=\pagetotal
        \advance\dimenone by \pagestretch
        \ifdim\dimenone<\pagegoal
          \dimenone=\pagetotal
          \advance\dimenone by#1         %addieren Skip vor Ueberschrift (=#1)
          \setbox0=\vbox{#2\parskip=0pt                %#2 ist gewaehlter Font
                     \hyphenpenalty=10000
                     \rightskip=0pt plus 5em
                     \noindent#3 \vskip#4}    %#3=Ueberschrift,#4=skip nachher
        \advance\dimenone by\ht0
        \advance\dimenone by 3\baselineskip   
        \ifdim\dimenone>\pagegoal\vfill\eject\fi
          \else\eject\fi\fi}

%%% OUR HEADLINE MACROS LOOK LIKE THIS USING THIS MACRO

\def\subheadline #1{\nin\bigbreak\vskip-\lastskip
      \checkleftspace{0.7cm}{\bf}{#1}{\medskipamount}
          \indent\vskip0.7cm\centerline{\bf #1}\medskip}

\def\sectionheadline #1{\bigbreak\vskip-\lastskip
      \checkleftspace{1.1cm}{\bf}{#1}{\bigskipamount}
         \vbox{\vskip1.1cm}\cen{\bfone #1}\bsk}

\def\lsectionheadline #1 #2{\bigbreak\vskip-\lastskip
      \checkleftspace{1.1cm}{\bf}{#1}{\bigskipamount}
         \vbox{\vskip1.1cm}\cen{\bfone #1}\msk \cen{\bfone #2}\bsk}

\def\lchapterheadline #1 #2{\bigbreak\vskip-\lastskip\indent\vskip3cm
                       \cen{\bftwo #1} \msk \cen{\bftwo #2} \gsk}
\def\llsectionheadline #1 #2 #3{\bigbreak\vskip-\lastskip\indent\vskip1.8cm
\cen{\bfone #1} \msk \cen{\bfone #2} \msk \cen{\bfone #3} \nobreak\bsk\nobreak}

%\def\[#1 #2\par{\hbox{\vtop{\hsize = 2.5 true cm \nin [#1]\hfill}
%\vtop{\hsize = 12.0 true cm \nin #2\penalty10000\llap.}}
%\vbox{\vskip.3\baselineskip}}  
% parameters for hsize are percentages !!! f.e. 0.2 + 0.8 = 1.0

\newtoks\literat
\def\[#1 #2\par{\literat={#2\unskip.}%
\hbox{\vtop{\hsize=.15\hsize\nin [#1]\hfill}
\vtop{\hsize=.82\hsize\nin\the\literat}}\par
\vskip.3\baselineskip}

\mathchardef\emptyset="001F 
\def\address{Author: \tt$\backslash$def$\backslash$address$\{$??$\}$}

\def\firstpage{\nin
{\obeylines \parindent 0pt }
\vskip2cm
\centerline {\bfone \title}
\gsk
\centerline{\bf\author}

\vskip1.5cm \rm}

\def\addresstwo{}

\def\dlastpage{\par\vbox{\vskip1cm\nin
\line{
\vtop{\hsize=.5\hsize{\parindent=0pt\baselineskip=10pt\nin\address}}
\quad 
\vtop{\hsize=.42\hsize\nin{\parindent=0pt
\baselineskip=10pt\addresstwo}}
\hfill} }}

% END OF LIEMACS.TEX

\pageno=1

\def\title{The $c$-function for non-compactly causal symmetric spaces}
\def\author{Bernhard Kr\"otz${}^*$ and Gestur \'Olafsson${}^{**}$}
\footnote{}{${}^*$Supported by the DFG-project HI 412/5-2

${}^{**}$Supported by LEQSF grant (1996-99)-RD-A-12}
\def\date{June 8, 1999}
\def\leftheadline{\tenbf\folio\hfil\eightrm\date}
\def\Box #1 { \msk\par\nin 
\centerline{
\vbox{\offinterlineskip
\hrule
\hbox{\vrule\strut\hskip1ex\hfil{\smc#1}\hfill\hskip1ex}
\hrule}\vrule}\msk }

\def\bs{\backslash} 

\def\address
{Bernhard Kr\"otz

Mathematical Institute 

TU Clausthal 

Erzstra\3e 1

D-38678 Clausthal--Zellerfeld

Germany

e-mail: mabk@math.tu-clausthal.de}

\def\addresstwo
{Gestur \'Olafsson

Department of Mathematics

Louisiana State University

Baton Rouge

LA 70803

e-mail: olafsson@math.lsu.edu}

\firstpage

\sectionheadline{Introduction} 
In this paper we prove a product formula for the $c$-function
associated to a non-compactly causal symmetric space ${\cal M}$.
Let us recall here the basic facts. Let $G$ be a connected semisimple
Lie group, $\tau :G\rightarrow G$ be a non-trivial involution and $H=G^{\tau}$.
Then ${\cal M}:=G/H$ is a semsimple symmetric space. The space ${\cal M}$
is called non-compactly causal, if $\q\:=\{X\in \g \: \tau (X)=-X\}$ contains
an open $H$-invariant hyperbolic cone $C\not= \emptyset$. In this case $S:=H\exp (C)$ is a
open subsemigroup of $G$. A {\it spherical function} on ${\cal M}$ is an $H$-biinvariant
continuous function on $S/H\subseteq {\cal M}$, which defines an eigendistribution
of the algebra of $H$-invariant differential operators on ${\cal M}$,
see [FH\'O94], [KN\'O98], [\'Ol97]. There exists a maximal abelian hyperbolic subspace
$\a$ of $\q$ such that $C=\Ad(H).(\a\cap C)$. Let 
$$\phi_{\lambda}(g.x_0)=\int_H a_H(gh)^{\lambda -
\rho}\, d\mu_H (h)$$
be a spherical function given by a convergent integral similar to the expression
for the spherical functions on the Riemannian symmetric spaces $G/K$. Here $x_0\in
{\cal M}$ is the coset $\{H\}$ and $a_H(g)\in A\:=\exp(\a)$ is determined by $g\in Ha_H(g)N$.    
The asypmtotic behaviour of $\phi_\lambda(a.x_0)$ along $S\cap A$ is given by
$\phi_\lambda (a.x_0)\sim c(\lambda )a^{\lambda-\rho}$, where $\rho$ is half
the sum over the positive roots counted with multiplicities. The function $c(\lambda)$
is the {\it $c$-function of the space ${\cal M}$}. It turns out that the
$c$-function is a product of two $c$-function, $c(\lambda )= c_{\Omega}(\lambda)c_0(\lambda)$
where $c_0(\lambda)$ is the Harish-Chandra $c$-function of a Riemannian subsymmetric
space $G(0)/K(0)$ and $c_{\Omega}(\lambda)$ is a function associated to
the real bounded symmetric domain $H/(H\cap K)$, where $K$ is a $\tau$-stable
maximal compact subgroup of $G$. The $c$-function was first introduced by
Oshima-Sekiguchi in [OS80], whereas $c_{\Omega}(\lambda)$ was first introduced
in [FH\'O94]. 

The $c$-function for a Riemannian symmetric space $G/K$ can be written as a
product of $c$-functions of rank one symmetric spaces associated 
to each  restricted root of $\g$ (Gindikin--Karpelevic formula).
For general non-Riemannian symmetric spaces $G/H$ one cannot expect this type of result.
However, for non-compactly causal symmetric spaces we show in this paper 
(cf.\ Theorem III.5) that such 
a product formula holds. The case of Cayley type spaces has already 
been treated by J.\ Faraut  
in [Fa95] by the use of Jordan algebra methods and in [Gr97] the case 
 $\Sl(n,\R)/\SO(p,q)$ is 
dealed. The approach presented here is general, different and 
relies on new insights on the fine convex geometry of the real bounded symmetric 
domain $\Omega$ (cf.\ Theorem II.5 and Theorem II.7.)

\par
Our result has important applications. The $c$-function was the last unknown part 
in the formula for the formal degree of the spherical holomorphic discrete series representations representations 
(cf.\ [Kr99]). Further it gives us  important information on 
the normalized spherical functions $\tilde \phi_\lambda \:= c_{\Omega}(\lambda)^{-1}
\phi_{\lambda}$. One knows that  
the function $\lambda \mapsto \tilde\phi_\lambda (s.x_0)$ has a meromorphic 
continuation to $\a_\C^*$ (cf.\ [\'Ol97]) and the product formula 
gives us important information on the poles. In particular, 
this allows more detailed analysis of the
spherical Laplace transform, in particular Paley-Wiener type theorems.

\sectionheadline{I. Non-compactly causal symmetric spaces and Lie algebras}

In this section we introduce notation and recall some facts concerning 
non-compactly causal symmetric Lie algebras and their associated 
symmetric spaces. Our source of reference is [Hi\'Ol96]. 

\subheadline{Algebraic preliminaries}

Let $\g$ be a simple finite dimensional real Lie algebra. Let $\tau : \g\rightarrow \g$ be
a non-trivial involution. Then 
$(\g,\tau)$  is a {\it symmetric Lie algebra}.  We write 
$\g=\h+\q$ for the $\tau$-eigenspace decomposition of $\g$ 
corresponding to the eigenvalues $+1$ and $-1$. Let $\theta$ be a
Cartan involution of $\g$ which commutes with $\tau$ and let 
$\g=\k+\p$ be the associated Cartan decomposition.

\par
For $\a,\b\subseteq \g$ let $\z_{\a}(\b)\:=\{X\in \a : [X,Y]=0, Y\in \b\}$ be the
{\it centralizer of $\b$ in $\a$}. We call $(\g,\tau)$ {\it non-compactly causal}, or
simply NCC, if 
$\z_{\q\cap\p}(\h\cap\k)\neq \{0\}$. We call $(\g,\tau)$ {\it non-compactly Riemannian
(NCR)} if $\tau$ is a Cartan involution.
If not otherwise stated from now on $(\g,\tau)$ denotes a
NCC symmetric Lie algebra.
Then $\z_{\q\cap\p}(\h\cap\k)=\z(\q\cap\p)=\R X_0$ is one dimensional. Let 
$\a\subeq \q\cap\p$ be a maximal abelian subspace and note that 
$\R X_0\subeq \a$ and that $\a$ is maximal abelian in $\p$. 
We write $\Delta=\Delta(\g,\a)$ for the root system of $\g$
with respect to $\a$ and 
$$\g=\z_\g(\a)\oplus\bigoplus_{\alpha\in \Delta} \g^\alpha$$
for the corresponding root space decomposition. 
We write 
$\g(0)\:=\h\cap\k+ \q\cap\p$ and note that 
$(\g(0), \tau(0))$, with $\tau(0)\:=\tau\res_{\g(0)}$, is 
NCR.
If $\alpha \in \Delta$ then either $\g^\alpha\subseteq \g(0)$
or $\g^{\alpha}\subseteq \q\cap \k + \h\cap\p$.
A root $\alpha\in\Delta$ is called {\it compact} if 
$\g^\alpha\subeq \g(0)$ and {\it non-compact} 
if $\g^\alpha\subeq \q\cap \k +\h\cap\p$. We write $\Delta_k$ 
and $\Delta_n$ for the collection of compact and non-compact 
roots, respectively. Note that $\Delta=\Delta_k\dot{\cup}\Delta_n$. 
\par
We can and will normalize $X_0$ such that $\Spec (\ad X_0)=\{-1,0,1\}$. 
Then $\Delta_k=\{\alpha\in\Delta\: 
\alpha(X_0)=0\}$ and we can choose a positive system 
$\Delta^+$ of $\Delta$ such that 
$$\Delta_n^+\:=\Delta_n\cap \Delta^+=\{\alpha\in \Delta_n\: 
\alpha(X_0)=1\}$$
and such that $\Delta_k^+\:=\Delta_k\cap \Delta^+$ is a positive
system in $\Delta_k$. Let $\Delta^-:=-\Delta^+$, $\Delta_n^-:=-\Delta_n^+$ and $\Delta_k^-:=-\Delta_k^+$.
\ssk  We recall now few facts about the structure of the root system $\Delta$.
Two roots $\alpha, \beta\in \Delta$ are said to be {\it strongly 
orthogonal} if $\alpha\pm \beta$ is not a root. 
Let $\Gamma\:=\{\gamma_1, \ldots,\gamma_r\}$ be a system
of strongly orthogonal roots in $\Delta_n^+$ of maximal length, 
i.e., $\Gamma$ consists of pairwise 
strongly orthogonal roots and has 
maximal number of elements with respect to this property.
We set 
$${\cal W}\:=N_{\Inn(\h\cap\k)}(\a)/ Z_{\Inn(\h\cap\k)}(\a)$$
and call ${\cal W}$ the {\it small Weyl group} of $\Delta$.

\Proposition I.1. For the root system $\Delta=\Delta(\g,\a)$ of a non-compactly 
causal symmetric Lie algebra $(\g,\tau)$ the following assertions hold:
\item{(i)} The root system $\Delta$ is reduced, 
i.e., if $\alpha\in\Delta$ then $2\alpha\not\in \Delta$. In particular, there
exists at most two root lengths. 
\item {(ii)} All long roots in $\Delta_n^+$ are conjugate under the small Weyl group
${\cal W}$. Moreover, all roots $\gamma_i$, $1\leq i\leq r$, are 
long.
\item{(iii)} Write $\Delta_{n,s}^+$ for the short roots in $\Delta_n^+$. Then, if 
$\Delta_{n,s}^+\neq \eset $, one has 
$$\Delta_{n,s}^+=\{ {1\over 2}(\gamma_i+\gamma_j)\: 1\leq i<j\leq r\}$$
and all elements of $\Delta_{n,s}^+$  are conjugate under ${\cal W}$.

\Proof. (i) [Hi\'Ol96, Th.\ 3.2.4] or [N\'O99, Lemma 2.12]. 
\par\nin (ii) [N\'O99, Lemma 2.26].
\par\nin (iii) [N\'O99, Lemma 2.22, Lemma 2.24].\qed 

For $\alpha \in \Delta$ let $H_{\alpha} \in \{[X,\tau(X)]\: X\in \g^\alpha\}\subeq \a$ 
be such that $\alpha (H_\alpha)=2$. For each $1\leq i\leq r$ let $H_i=H_{\gamma_i}$.
We set $\cc\:=\span_\R\{H_1,\ldots, H_r\}\subeq \a$ and write $\b$ for the 
orthogonal complement of $\cc$ in $\a$ with respect to the Cartan-Killing form; 
in particular $\a=\cc\oplus\b$.

\Proposition I.2. The positive system $\Delta_k^+$ can be chsosen 
such that for the restriction of $\Delta=\Delta(\g,\a)$ to $\cc$ the 
following assertions hold:
$$\Delta_n^+\res_\cc=\{ {1\over 2}(\gamma_i+\gamma_j)\: 1\leq i,j\leq r\}\cup\{{1\over 2}\gamma_i\: 
1\leq i\leq r\},$$
$$\Delta_k^+\res_\cc\bs\{0\}=
\{ {1\over 2}(\gamma_i-\gamma_j)\: 1\leq j<i \leq r\}
\cup\{-{1\over 2}\gamma_i\: 
1\leq i\leq r\}.$$
Moreover, the second sets in the two unions from above may or may not  occur simultaneously.

\Proof. [N\'O99, Th.\ 2.21] or [Kr99, Th.\ IV.4].\qed

Since we have free choice for $\Delta_k^+$ we assume in the sequel that 
$\Delta_k^+\res_\cc\bs\{0\}\subeq\{ {1\over 2}(\gamma_i-\gamma_j)\: 1\leq j<i \leq r\}\cup \{-{1\over 2}\gamma_i\: 1\leq i\leq r\}$.

\Lemma I.3. Assume that $\Delta_{n,s}\neq\eset$ and let $\Pi_k$ be the 
set of simple roots corresponding to $\Delta_k^+$. Then there exits $\beta_1,\ldots,\beta_m\in\b^*$,
$\beta_j(X_0)=0$, and $\delta_1,\ldots\delta_l\in\b^*$, $\delta_i(X_0)={1\over 2}$, such that 
$$\Pi_k=\{ {1\over 2} (\gamma_{i+1}-\gamma_i)\: 1\leq i\leq r-1\}\cup\{\beta_1,
\ldots,\beta_m\} \cup \{-{1\over 2}\gamma_r +\delta_i\: 1\leq i\leq l\}.$$
Here the last set occurs if and only if there exits half roots in 
$\Delta\res_\cc$.

\Proof. For each $\alpha\in\Delta$ let $s_\alpha$ denote the corresponding 
reflection. Then $s_{\gamma_j}({1\over 2}(\gamma_i+\gamma_j))={1\over 2}(\gamma_i-\gamma_j)$,
$i\neq j$ together with Proposition I.1(iii) shows that 
$\Delta_k\supeq\{ {1\over 2}(\gamma_i-\gamma_j)\: 1\leq i\neq j\leq r\}$. 
Thus Proposition I.2 yields that 
$$\eqalign {\Delta_k^+ &\subeq \{{1\over 2}(\gamma_i-\gamma_j)\: 1\leq j<i \leq r\}+
(\b^*\cap X_0^\bot)\cr 
&\cup \{-{1\over 2}\gamma_i\: 1\leq i\leq r\}+\{\delta
\in\b^*\: \delta(X_0)={1\over 2}\}\cup \b^*.\cr}$$
Now  the assertion follows easily from 
Proposition I.2 and the fact that $\Delta$ is a root system. \qed

We define the {\it maximal cone} in $\a$ is defined by 
$$C_{\rm max}\:=\{X\in \a\: (\forall \alpha\in\Delta_n^+)\,
\alpha(X)\geq 0\}.$$

\Lemma I.4. Let $X_0=X_0^b+X_0^c$ with $X_0^b\in \b$ and $X_0^c\in\cc$. Then we have 
$X_0^b, X_0^c\in C_{\rm max}$. 

\Proof. First note that $X_0^c={1\over 2}(H_1+\ldots+H_r)$ and so $X_0^c\in C_{\rm max}$ by 
Proposition I.2. 
\par\nin To show $X_0^b\in C_{\rm max}$ let $\alpha\in \Delta_n^+$. Then Proposition I.2
shows that $\alpha={1\over 2}(\gamma_i+\gamma_j)+\beta$ with $\beta\in\b^*$, $\beta(X_0)=\beta
(X_0^b)=0$, or $\alpha=-{1\over 2}\gamma_i +\delta$ with $\delta\in \b^*$ and $\delta(X_0)
=\delta(X_0^b)={1\over 2}$. In any case we have $\alpha(X_0^b)\geq 0$ concluding 
the proof of the lemma.\qed

Finally we define subalgebras of $\g$ by 
$$\n\:=\bigoplus_{\alpha\in\Delta^+} \g^\alpha, \quad \oline \n\:=
\bigoplus_{\alpha\in \Delta^-} \g^\alpha, \quad \n_k^\pm\:=
\bigoplus_{\alpha\in\Delta_k^\pm}\g^\alpha, \quad \n_n^\pm\:=\bigoplus_{\alpha\in 
\Delta_n^\pm} \g^\alpha$$ 
and note that $\n=\n_n^+\rtimes \n_k^+$ and $\oline \n=\n_n^-\rtimes \n_k^-$
are semidirect products.

\subheadline{Analytic preliminaries}

Let $G_{\C}$ be a simply connected Lie group with Lie algebra $\g_{\C}$ and
let $G$ be the analytic subgroup of $G_{\C}$ corresponding to $\g$. Let
$H=G^\tau=\{X\in G\: \tau (g)=g\}$. We write 
$A$, $K$, $N$, $\oline N$, $N_k^\pm$, $N_n^\pm$ for the analytic
subgroups of $G$ which correspond to $\a$, $\g(0)$, $\h$, $\k$, $\n$, 
$\oline \n$, $\n_k^\pm$, $\n_n^\pm$. Note that the groups $A$, $N$, 
$\oline N$, $N_k^\pm$, $N_n^\pm$ are all simply connected and that the 
corresponding exponential mappings $\exp_A\: \a\to A$, $\exp_N\: \n\to N$ 
etc. are all diffeomorphisms. Let $G(0) = Z_G(X_0)=
\{g\in G\: \Ad (g).X_0=X_0\}$. Then $H$ and $G(0)$ are $\tau$ and $\theta$
invariant, $H=(H\cap K) \exp(\h\cap \p)$ and $G(0)= (H\cap K) \exp (\q\cap \p)$.

\par The Lie algebra  $\g$ decomposes as $\g=\h+\a+\n$ and the multiplication 
mapping 
$$H\times A\times N\to G, \ \ (h,a,n)\mapsto han$$
is an analytic diffeomorphism onto its open image $HAN$. 

\par Note that $\oline N= N_n^-\rtimes N_k^-$. We have 

$$\oline N\cap HAN=\exp(\Omega) N_k^-=N_k^-\exp(\Omega)\leqno(1.1)$$
with $\Omega\cong H/ H\cap K$ a real bounded symmetric domain in $\n_n^-$.

\sectionheadline{II. The geometry of the real bounded symmetric domain $\Omega$}

We denote by $\kappa$ the Cartan-Killing form on $\g$ and define an 
inner product on $\g$ by $\la X,Y\ra \:=-\kappa(X,\theta (Y))$ for $X,Y\in \g$.
Let $X_i\in \g^{\gamma_i}$ be such that $H_i=[X_i,X_{-i}]$, with $X_{-i}=\tau
(X_i)$.
By [Hi\'Ol96] and Herman's Convexity Theorem we have
$$ \leqalignno{\Omega &=\{X\in \n_n^-\: \|\ad (X+\tau(X))\|<1\}&(2.1)\cr
&= \Ad (H\cap K).\{\sum_{j=1}^rt_jX_{-j}\: -1<t_j<1, j=1,\ldots ,r\}\, ,&(2.2)
\cr}$$
where $\|\cdot\|$ denotes the 
operator norm corresponding to the scalar product $\la\cdot,\cdot\ra$ on $\g$.
Note that (2.1) implies that $\Omega$ is a convex balanced subset of $\n_n^-$.

\Remark II.1. Recall the definition of the  maximal cone $C_{rm max}$  in $\a$.
Then it is clear from the characterization (2.1) of $\Omega$ that 
$e^{\ad X}.\Omega\subeq \Omega$ for all $X\in C_{\rm max}$. We also 
have a {\it minimal cone} in $\a$ defined by  
$$C_{\rm min}\:=\cone(\{ [X,\tau(X)]\: X\in\g^\alpha, \alpha\in \Delta^+\})
=\overline{\sum_{\alpha\in  \Delta_n^+}\R^+ H_{\alpha} }\, 
.$$
We note that $C_{\rm min}\subeq C_{\rm max}$ and  
in particular $H_i\in C_{\rm max}$ for each $1\leq i\leq r$.\qed

The following concept turns out to be very useful for the investigation of the 
fine convex geomety of $\Omega$.

\Definition II.2. (Oshima-Sekiguchi) By a {\it signature  of} $\Delta$
we understand a map $\eps\: \Delta\to \{-1,1\}$ with the following 
properties:
\item{(S1)} $\eps(\alpha)=\eps(-\alpha)$ for all $\alpha\in \Delta$.
\item{(S2)} $\eps(\alpha+\beta)=\eps(\alpha)\eps(\beta)$ for all 
$\alpha,\beta\in\Delta$ with $\alpha+\beta\in
\Delta$.\qed
\par
If $\eps : \Delta \to \{-1,1\}$ is a signature then $\theta_\eps : \g\to \g$ defined by
$\theta_\epsilon (X)=\eps (\alpha)\theta (X)$, $X\in \g^\alpha$ and
$\theta_\epsilon|_{\z_{\g}(\a)}=
\theta|_{\z_{\g}(\a)}$ is an involution on $\g$ that commutes with $\theta$
(see [OS80, Def.\ 1.2]). As $\tau|_{\z_{\g}(\a)}=\theta|_{\z_{\g}(\a)}$ and
$\tau|_{\g^\alpha}=\pm \theta|_{\g^\alpha}$, with $+$ if $\alpha$ is compact and
$-1$ if $\alpha$ non-compact, it follows that $\theta_\eps$ also commutes with
$\tau$.

\Lemma II.3. Keep the notation of Definition II.2. 
\item{(i)} If $\eps$ is a signature of $\Delta$, then 
the prescription 
$$\sigma_\eps(X)\:=\cases{X & for $X\in\z_\g(\a)$,\cr \eps(\alpha)X & for $X\in 
\g^\alpha$, $\alpha\in\Delta$\cr}$$
defines an involutive automorphism of $\g$. The involution $\sigma_\epsilon$ commutes with
both $\tau$ and $\theta$.
\item{(ii)} Let
$\Pi\:=\{\alpha_1, \ldots,\alpha_n\}$ be a basis of $\Delta$. Then 
for any collection $(\eps_1,\ldots, \eps_n)\in\{-1,1\}^n$ one can define 
a signature $\eps$ of $\Delta$ by setting 
$$\eps(\pm\sum_{i=1}^n n_i\alpha_i)\:=\prod_{i=1}^n \eps_i^{n_i}\quad
\hbox{for} \quad \sum_{i=1} n_i\alpha_i \in\Delta.$$
\item{(iii)} Let the notation be as in (ii). Then
$\eps \mapsto (\eps (\alpha_i))_{i=1}^n$ defines a bijection between the set of signatures of $\Delta$ and 
$\{-1,1\}^n$.

\Proof. (i) This follows by the Oshima-Sekiguchi construction because $\sigma_\epsilon
=\tau_\epsilon\theta$. (ii) is clear and (iii) follows from (ii).  \qed 
\par
In the sequel we identify signatures with elements in 
$\{-1,1\}^n$.

\Lemma II.4. Let $\eps$ be a signature of $\Delta$. Then $\sigma_\eps(\Omega)=
\Omega$.

\Proof.  Let $X\in \Omega$. By (2.2) there is a $k\in H\cap K$ and
$Y=\sum_{j=1}^rt_jX_{-j}$, $-1<t_j<1$ such that
$\Ad(k).Y=X$. As $\sigma_\eps$ commutes with $\tau$ and $\theta$
it follows that $\sigma_\eps (k)\in K\cap H$. Hence
$\sigma_\eps (X)=\Ad (\sigma_\eps(k)).\sum_{j=1}^r\eps (\gamma_j)t_jX_{-j}\in
\Omega$.\qed 

Recall that there is basis $\Pi\subeq \Delta^+$ having the form 
$$\Pi=\{\alpha_0, \alpha_1, \ldots, \alpha_n\}$$
with $\alpha_0$ long and non-compact and $\alpha_i$, $1\leq i\leq n$ compact. Thus every 
non-compact negative root $\gamma\in \Delta_n^-$
can be written as $\gamma=-\alpha_0-\sum_{i=1}^n m_i\alpha_i$, 
$m_i\in\N_0$. By our choice of $\Delta_k^+$ we have $\alpha_0=\gamma_1$.

\Theorem II.5. For each $\gamma\in\Delta_n^-$ let $p_\gamma\: \n_n^-\to \g^\gamma$
be the orthogonal projection. Then 
$$X\in\Omega\Rarrow p_\gamma(X)\in\Omega.$$

\Proof. Let $X=\sum_{\gamma\in \Delta_n^-} X_\gamma\in\Omega$ with 
$X_\gamma\in \g^\gamma$, $\gamma\in \Delta_n^-$.
We have  to show that $X_\gamma\in \Omega$. 
Recall that there are at most two root length in $\Delta$ (cf.\ Proposition 
I.1(i)).
\par\nin Case 1: $\gamma$ is a long root.
\par By Proposition I.1(ii) there exists an element $h\in N_{\Inn(\h\cap\k)}
(\a)$ such that $h.\gamma=-\alpha_0$. Thus we may assume that $\gamma=-\alpha_0=-\gamma_1$.
Let $H\:=\sum_{j=2}^r H_j$. By Remark II.1 we
have 
$$X_1\:=\lim_{t\to+\infty} e^{t\ad H}.X\in \Omega.$$
If we express $X_1=\sum_{\beta\in\Delta_n^-} X_\beta$ 
as a sum of root vectors, then Proposition I.2 implies that 
$\beta\res_\cc=-\gamma_1$ or $\beta={1\over 2}\gamma_1-\delta$ with 
$\delta(X_0^b)={1\over 2}$.
Since $X_0^b\in C_{\rm max}$ (cf.\ Lemma I.4), we now get 
$$X_\gamma=\lim_{t\to+\infty} e^{t\ad X_0^b}.X_1\in \Omega.$$ 

\par\nin Case 2: $\gamma$ is a short root. 
\par By Proposition I.1(iii) we may assume that $\gamma={1\over 2}(\gamma_1+\gamma_2)$ and by 
Lemma I.3 we may suppose $\alpha_0=\gamma_1$, $\alpha_j={1\over 2}(\gamma_{j+1}-\gamma_j)$
for $1\leq j\leq r-1$. 
Write   
$$X=\sum_{m_i\geq 0} X_{m_1,\ldots,m_n},$$
where $X_{m_1,\ldots, m_n}\in \g^{-(\alpha_0+\sum_{i=1}^n
m_i\alpha_i)}$. Then we  have to show that $X_{1,0,\ldots,0}\in \Omega$.
Set 
$$X_{\rm ev}\:=\sum_{ m_n\equiv 0(2)} X_{m_1,\ldots, m_n}\quad\hbox{and}
\quad X_{\rm odd}\:=\sum_{m_n\equiv 1(2)} X_{m_1,\ldots, m_n}.$$
Then $X=X_{\rm ev} +X_{\rm odd}$ and we claim that 
$X_{\rm ev}$, $X_{\rm odd}\in 
\Omega$.  Let $\eps=(1,1,1,\ldots, -1)$. Then by Lemma II.4 we get:
$$\sigma_\eps(X)=\sigma_\eps(X_{\rm ev} +X_{\rm odd})=X_{\rm ev} -X_{\rm odd}
\in \Omega\, .$$ 
Since $\Omega$ is balanced and convex we moreover have 
$$X_{\ev}={1\over 2}(X+\sigma_{\eps}(X))\in\Omega\quad\hbox{and}\quad
X_{\rm odd}={1\over 2}(X-\sigma_{\eps}(X))\in\Omega.$$
By repeating this argument we thus my assume that 
$$X=\sum_{m_1\equiv 1(2)\atop m_j\equiv 0(2), j>1}  X_{m_1,\ldots,m_n}.$$
Now we apply the contraction semigroup generated by  
$H=\sum_{j=3}^r H_j\in C_{\rm max}$ and obtain
$$X_1\:=\lim_{t\to+\infty} e^{t\ad H}.X\in \Omega.$$
Thus we may assume $X=X_1$ and $X=\sum_{\beta\in\Delta_n^-} X_\beta$  
with $-\beta=\gamma, \gamma_1,\gamma_2,
{1\over 2}(\gamma_1+\gamma_2)+\beta, -{1\over 2}\gamma_1+\sigma_1, -{1\over 2}\gamma_2+\sigma_2$
and $\beta,\sigma_1, \sigma_2
\in\b^*$, $\sigma_1(X_0)=\sigma_2(X_0)={1\over 2}$ (cf.\ Proposition I.2).
Write $\beta=-\gamma_1-\sum_{j=1}m_j\alpha_j$. 
The cases $\beta=\gamma_1$ and $\beta=\gamma_2$ are excluded,
since we have $m_1=0$, resp. 
$m_1=2$, contradicting $m_1\equiv 1(2)$. Applying to $X$ the 
contraction semigroup generated by $X_0^b\in C_{\rm max}$ excludes the case 
$\beta=-{1\over 2}\gamma_1+\sigma_1$ and $\beta =-{1\over 2}\gamma_2+\sigma_2$. 
Let now $Y\in\b$ such that $\delta_j(Y)>0$, $1\leq j\leq l$, and 
$\beta_j(Y)>0$, $1\leq j\leq m$ (cf.\ Lemma I.3). 
Then $\Delta_n^+\subeq \N_0[\Pi]$ shows that 
$Y\in C_{\rm max}$. But then 
$$X_\gamma=\lim_{t\to+\infty} e^{t\ad Y}.X\in \Omega,$$
completing the proof Case 2 and hence of the theorem. \qed

\subheadline{Subdomains of rank one}

For $\alpha\in\Delta^+$ we set 
$$\g(\alpha)\:=\big(\g^\alpha+\g^{-\alpha}+[\g^\alpha, \g^{-\alpha}]\big)'$$
and  $\tau(\alpha)\:=\tau\res_{\g(\alpha)}$. Then $(\g(\alpha),
\tau(\alpha))$ is a symmetric  subalgebra of $(\g,\tau)$ of real 
rank one, that is $\a(\alpha)\:=\a\cap\g(\alpha)$ is one dimensional.
Further we set $\h(\alpha)\:=\h\cap \g(\alpha)$ etc. 
We denote by $G(\alpha)$, $A(\alpha)$
etc. the analytic subgroups of $G$ corresponding to $\g(\alpha)$, $\a(\alpha)$ etc.
Let $H(\alpha )=G(\tau)^{\tau (\alpha)}=G(\alpha)\cap H$.
\par Assume that $\alpha\in\Delta_n^+$. Then $(\g(\alpha), \tau(\alpha))$ is 
NCC and 
$\n(\alpha)=\n_n^+(\alpha)=\g^\alpha$. Let $\Omega(\alpha)\cong H(\alpha)/
\big (K(\alpha)\cap H(\alpha)\big)$ 
be the real bounded symmetric domain in $\oline\n(\alpha)=\n_n^-(\alpha)$.

\Lemma II.6. Let $\alpha\in \Delta_n^+$ and $s_\alpha\in G(\alpha)$ be a representaive of the one element 
big Weyl group $N_{G(\alpha)}(\a(\alpha))/ Z_{G(\alpha)}(\a(\alpha))$ of $\g(\alpha)$. Then 
$$\big(\oline N(\alpha)\cap H(\alpha) A(\alpha) N(\alpha)\big)
\dot{\cup}\big(\oline N(\alpha)\cap H(\alpha) s_\alpha A(\alpha) N(\alpha)\big 
)$$
is open and dense in $\oline N(\alpha)$. 

\Proof. This follows by
Matsukis Theorem (cf.\ [Ma79, Theorem 3]), if we can show that
$M(\alpha)\:=Z_{K(\alpha )}(\a (\alpha))\subseteq H(\alpha)$
because $s_\alpha M(\alpha)
=M(\alpha )s_\alpha$. 
Let $F=\exp(i\a (\alpha ))\cap G(\alpha)$. Then one has 
$M(\alpha )=FZ_{H(\alpha)_o}(\a (\alpha))$ by [N\'O99,\ Lemma 5.7]. But if
$f\in F$ then $\tau(\alpha)(f)=f^{-1}=f$, by the same lemma. Hence
$F\subseteq H(\alpha)$, which implies that $M(\alpha)\subseteq H(\alpha)$.
\qed

\Theorem II.7. Let $\alpha\in\Delta_n^+$. Then 
$\Omega\cap \oline\n (\alpha)=\Omega(\alpha)$.

\Proof. "$\supeq$": This is clear.
\par\nin "$\subeq$": Note that $\Omega\cap\oline\n(\alpha)$ is open and convex in 
$\oline\n(\alpha)$. We have  
$$\exp(\Omega)\cap \big(H(\alpha) s_\alpha A(\alpha) N(\alpha)\big)=\eset, \leqno(2.3)$$
since $\exp(\Omega)\subeq HAN$ and $HAN\cap H s_\alpha AN=\eset $ by 
Matsukis Theorem. In view of (2.3), Lemma II.7 implies that 
there exists an open dense subset $\Omega_\alpha$ of $\Omega\cap\oline \n(\alpha)$ such 
that $\Omega_\alpha\subeq \Omega (\alpha)$. Now the assertion follows from the fact 
that both $\Omega(\alpha)$ and $\Omega\cap \oline \n(\alpha)$ are open and convex.
\qed

\sectionheadline {III. The product formula for the $c$-function}

Recall the $HAN$-decomposition in $G$ from Section I. For each $\lambda\in\a_\C^*$
and $g$ in $G$ we set
$$a_H(g)^\lambda \:=\cases{0 & if $g\not\in HAN$,\cr e^{\lambda(\log a)} & if 
$g=han \in HAN$}.$$

\par For a locally compact group $G$ we write $\mu_G$ for a left Haar 
measurwe on $G$. 

\Definition III.1. (The $c$-functions) For each $\alpha\in\a^*$ let 
$m_\alpha\:=\dim \g^\alpha$ and put $\rho\:={1\over 2}\sum_{\alpha\in \Delta^+} m_\alpha
\alpha$. For $\lambda\in\a_\C^*$ we now set 

$$c(\lambda)\:=\int_{\oline N} a_H(\oline n)^{-(\lambda+\rho)}\ d\mu_{\oline N}(\oline n)
=\int_{\oline N\cap HAN} a_H(\oline n)^{-(\lambda+\rho)}\ d\mu_{\oline N}(\oline n),$$
$$c_{\Omega}(\lambda)\:=\int_{N_n^-} a_H(\oline n)^{-(\lambda+\rho)}
\ d\mu_{N_n^-}(\oline n)=\int_{\Omega} a_H(\oline n)^{-(\lambda+\rho)}
\ d\mu_{N_n^-}(\oline n), $$
and 
$$c_0(\lambda)\:=\int_{N_k^-} a_H(\oline n)^{-(\lambda+\rho)}\ d\mu_{N_k^-}(\oline n)$$
whenever the defining integral exist. 
We write ${\cal E}$, ${\cal E}_\Omega$ and 
${\cal E}_0$ for the domain of definition of $c$, $c_\Omega$ and $c_0$, respectively.
We call $c$ the $c${\it -function of the non-compactly causal symmetric space} $G/H$ and
$c_\Omega$ the $c${\it -function of the real bounded symmetric domain} $\Omega$, 
while $c_0$ is the usual $c$-{\it function of the non-compact Riemannian
symmetric space} $G(0)/K(0)$.\qed

\Remark III.2. (a) The choice of the particular analytic realization $G/H$ of 
$(\g,\tau)$ as a symmetric space is immaterial for the definition of the 
$c$-function. 
\par\nin (b) We have ${\cal E}={\cal E}_0\cap{\cal E}_\Omega$ and for all
$\lambda\in {\cal E}$ one has the splitting 
$$c(\lambda)=c_0(\lambda)c_\Omega(\lambda)$$   
(cf.\ [FH\'O94, Lemma 9.2]).
\par\nin (c) The $c$-functions can be written as Laplace transforms (cf.\ [KN\'O98]). 
Let us explain this for the $c$-function $c$. For $c_0$ and $c_\Omega$ one has analogous
statements. 
\par  There exists a positive Radon measure $\mu$ on $\a$ such that 
$$(\forall \lambda\in {\cal E})\quad  c(\lambda)={\cal L}_\mu(\lambda)\:=
\int_\a e^{\lambda(X)} \ d\mu(X), $$
i.e., $c$ is the Laplace transform of $\mu$. In particular we see that 
the domain of definition ${\cal E}$ is a tube 
domain over a convex set, i.e., one has
$${\cal E}=i\a^*+ {\cal E}_\R$$
with ${\cal E}_\R\subeq \a^*$ a convex subset of $\a^*$. One knows that 
$\Int {\cal E}$ is non-empty. Moreover, the fact that $c$ is a Laplace transform 
implies that 
$c$ is holomorphic on $\Int {\cal E}$ and that $c$ has no holomorphic extension to a
connected open tube domain strictly larger than $\Int {\cal E}$.\qed

Now we are going to prove the product formula for the $c$-function $c_\Omega$. 
Our srategy is a modified Gindikin-Karpelevic approach as presented in 
[GaVa88, p.\ 175--177] or [Hel84, Ch.\ IV].

For a positive system $R\subeq \Delta$ we set $\oline \n_R\:=
\bigoplus_{\alpha\in -(\Delta^+\cap R)} \g^\alpha$ and write 
$\oline N_R$ for the corresponding analytic subgroup of $G$. We define 
an auxiliary $c$-function by 

$$c_R(\lambda)\:=\int_{\oline N_R} a_H(\oline n)^{-(\lambda+\rho)} \ d\mu_{\oline N_R} (\oline n)$$
whenever the integral exists. 
 
\par For a single root $\alpha\in\Delta^+$ we set $\rho_\alpha\:={1\over 2}m_\alpha\alpha$ 
and write 

$$c_\alpha(\lambda)\:=\int_{\oline N(\alpha)} a_{H(\alpha)}(\oline n)^{-(\lambda+\rho_\alpha)}
\ d\mu_{\oline N(\alpha)} (\oline n).$$ 
We denote by ${\cal E}_\alpha\subeq \a_\C^*$ the domain of definition 
of $c_\alpha$.

\Proposition III.3. For any positive system $R\subeq \Delta$ we have that 
$$c_R(\lambda)=\prod_{\alpha\in (R\cap\Delta^+)} c_\alpha(\lambda)$$
and $c_R(\lambda)$ is defined if and only if $\lambda\in\bigcap_{\alpha\in 
(R\cap\Delta^+)} {\cal E}_\alpha$.

\Proof. We proceed by induction on $|R\cap\Delta^+|$. If $R\cap \Delta^+=\eset$, 
then the assertion is clear. 
\par Assume that $R\cap \Delta^+\neq \eset$. Then we find an element 
$\beta\in R\cap\Delta^+$ which is simple in $R$. Set 
$Q\:=s_\beta.R$. Then $Q=R\bs \{\beta\} \cup \{ -\beta\}$ since 
$\Delta$ is reduced (cf.\ Proposition I.1(i)). Thus we have 
$(Q\cap\Delta^+)\dot\cup\{\beta\}=R\cap\Delta^+$.
We now have to distinguish to cases. 

\par\nin Case 1: $\beta$ is compact. 
\par\nin In this case, the $HAN$-decomposition of $G(\beta)$ coincides with 
the Iwasasa decomposition, i.e.  $G(\beta)=K(\beta)A(\beta)N(\beta)$. Thus
$c_R(\lambda)=c_\beta(\lambda) c_Q(\lambda)$ follow as in [GaVa88, Prop.\ 4.7.6].

\par\nin Case 2: $\beta$ is non-compact. 
\par\nin Set $\oline N_Q^k\:=\oline N_Q\cap N_k^-$, $\oline N_Q^n\:=
\oline N_Q\cap N_n^-$ and note that $\oline N_Q=\oline N_Q^n\rtimes \oline  N_Q^k$.
Since $\oline N_R=\oline N(\beta) \oline N_Q$ we thus get 

$$c_R(\lambda)=\int_{\oline N(\beta)} \int_{\oline N_Q^n} \int_{\oline N_Q^k} 
a_H(\oline n_\beta\oline n_n \oline n_k)^{-(\lambda+\rho)} 
\ d\mu_{\oline N(\beta)} (\oline n_\beta)\ d\mu_{\oline N_Q^n} (\oline n_n)
\ d\mu_{\oline N_Q^k }(\oline n_k).$$

If $\oline n_\beta \oline n_n \oline n_k\in \oline N\cap HAN$, then 
(1.1) implies that $\oline n_\beta\oline n_n\in \exp(\Omega)$. Since 
$\n_n^-$ is abelian, Theorem II.5 therefore implies that $\oline n_\beta\in \exp(\Omega)$
and so $\oline n_\beta\in \exp(\Omega(\beta))$ by Theorem II.7. Therefore 
we can write $\oline n_\beta =h_\beta a_\beta n_\beta$ with $h_\beta\in H(\beta)$, 
$a_\beta\in A(\beta)$ and $n_\beta\in N(\beta)$. Now one can proceed as 
in [GaVa88, p.\ 175--177] and one gets 
$c_R(\lambda)=c_\beta(\lambda) c_Q(\lambda)$.\qed

\Remark III.4. If we choose $R=-\Delta_n^+\cup \Delta_k^+$ (this is a positive 
system since $\Delta_n^+$ is ${\cal W}$-invariant), then we have 
$c_0=c_R$ and Proposition III.3 results in the Gindikin-Karpelevic product 
formula 
$$c_0(\lambda)=\prod_{\alpha\in\Delta_k^+} c_\alpha(\lambda)$$
of the $c$-function $c_0$ on $G(0)/ K(0)$ (cf.\ [GaVa88, Th.\ 4.7.5]
or [Hel84, Ch.\ IV, Th.\ 6.13, 6.14]).\qed

\Theorem III.5. {\rm (The product formula for $c_\Omega$)} For the $c$-function 
$c_\Omega$ of the real bounded symmetric domain $\Omega$ one has
$${\cal E}_\Omega=\{ \lambda\in \a_\C^*\: (\forall \alpha\in \Delta_n^+)
\Re \lambda(H_\alpha)< 2-m_\alpha\}$$
and 
$$c_\Omega(\lambda)=\kappa\prod_{\alpha\in\Delta_n^+} B\Big({m_\alpha\over 2}, 
-{\lambda(H_\alpha)\over 2} -{m_\alpha\over 2} +1\Big) $$
where $B$ denotes the Beta function and $\kappa$ is a positive constant only 
depending on $(\g,\tau)$.

\Proof. Set ${\cal E}_\Omega'\:=\bigcap_{\alpha\in\Delta_n^+}{\cal E}_\alpha$. 
We want to apply Proposition III.3 to $R=\Delta^+$. In view of Remark III.2(b) and 
Remark III.4, we thus  get 
$$(\forall \lambda\in {\cal E}\cap{\cal E}_\Omega' )\quad c_\Omega(\lambda)=
\prod_{\alpha\in\Delta_n^+} c_\alpha(\lambda)=
\prod_{\alpha\in\Delta_n^+} c_{\Omega(\alpha)}(\lambda).\leqno(3.1)$$
By [FH\'O94, (10.3)] one has 
$$c_{\Omega(\alpha)}(\lambda)= 2^{m_\alpha-1} B\Big({m_\alpha\over 2}, 
-{\lambda(H_\alpha)\over 2} -{m_\alpha\over 2} +1\Big)\leqno(3.2)$$
and 
$${\cal E}_{\Omega(\alpha)}=\{\lambda\in\a_\C^*\: \Re\lambda(H_\alpha)< 2-m_\alpha\}
.\leqno(3.3)$$
It follows from (3.3) that 
$${\cal E}_\Omega'=\{ \lambda\in \a_\C^*\: (\forall \alpha\in \Delta_n^+)
\Re \lambda(H_\alpha)< 2-m_\alpha\}.\leqno (3.4)$$
Besides ${\cal E}_\Omega={\cal E}_\Omega'$ all assertions of the theorem now 
follow from (3.1)-(3.4). 
Finally, ${\cal E}_\Omega= {\cal E}_\Omega'$ follows from the fact that all
$c$-functions involved are Laplace transforms (cf.\ Remark III.2(c)).\qed 

\par\nin 
The following simple fact that shows that we can split off all the non-compact
roots to get the $c_{\Omega}$-function before we come to the compact roots.

\Lemma III.6. Let $R$ be a any positive system of roots in $\Delta$. If 
$R\cap \Delta_n^+\neq \eset$, then $R\cap \Delta_n^+$ contains a 
a root that is simple in $R$.

\Proof. Let $\{\beta_0,\ldots, \beta_n\}$ be the set of simple roots 
in $R$. Let $\gamma\in R\cap\Delta_n^+$. Then 
$\gamma=\sum_{i=0}^n n_i \beta_i$ with $n_i\in \N_0$. Thus 
$1=\gamma(X_0)=\sum_{i=0}^n n_i \beta_i(X_0)$ which implies that $\beta_i(X_0)>0$
for at least one $\beta_i$. But then $\beta_i\in \Delta_n^+$.\qed

\def\entries{

\[Fa95 Faraut, J., {\it Fonctions Sph\'eriques sur un Espace Sym\'etrique
Ordonn\'e de Type Cayley}, Contemp. Math. {\bf 191} (1995), 41--55

\[FH\'O94 Faraut, J., J.\ Hilgert, and G. \'Olafsson, {\it Spherical functions
on ordered symmetric spaces}, Ann. Inst. Fourier {\bf 44} (1994), 927--966

\[GaVa88 Gangolli, R., and V.S.\ Varadarajan, ``Harmonic Analysis of 
Spherical Functions on Real Reductive Groups,'' Ergebniss der
Mathematik {\bf 101}, Springer, 1988 

\[Gr97 Graczyk, P., {\it Function $c$ on an ordered symmetric space}, Bull. Sci.
math. {\bf 121} (1997), 561--572

\[Hel84 Helgason, S., "Groups and Geometric Analysis", Acad. Press, London, 1984

\[Hi\'Ol96 Hilgert, J.\ and 
G.\ \'Olafsson, ``Causal Symmetric Spaces, Geometry and
Harmonic Analysis,'' Acad. Press, 1996 

\[Kr99 Kr\"otz, B., {\it Formal dimension of semisimple symmetric spaces},
Compositio math., to appear 

\[KN\'O98 Kr\"otz, B., K.--H. Neeb, and G.\ \'Olafsson, {\it Spherical Functions
on Mixed Symmetric Spaces}, submitted 

\[Ma79 Matsuki, T., {\it The orbits of affine symmetric spaces under the action 
of minimal parabolic subgroups}, J. Math. Soc. Jpn. {\bf 31}, 331--357 (1979)

\[N\'O99 Neumann, A., and G. \ \'Olafsson, 
{\it Minimal and Maximal Semigroups Related to Causal Symmetric 
Spaces}, Semigroup Forum, to appear 

\[\'Ol97 \'Olafsson, G., {\it Spherical Functions and
Spherical Laplace Transform on Ordered Symmetric Spaces}, submitted

\[OS80 Oshima, S., Sekiguchi, J, {\it Eigenspaces of Invariant Differential
Operators on an Affine Symmetric Spaces}, Invent. math. {\bf 57} (1980), 1--81
}
{\sectionheadline{\bf References}
\frenchspacing
\entries\par}
\dlastpage
\hfill\eject
\bye